\documentclass[letterpaper, 10 pt, conference]{ieeeconf}  

\usepackage{graphics} 
\usepackage{epsfig} 
\usepackage{mathptmx} 
\usepackage{amsfonts}
\usepackage{amsmath} 
\usepackage{amssymb}  
\usepackage{stdmath}  
\usepackage{subfigure}
\usepackage{multirow}
\usepackage{xcolor,colortbl}
\usepackage[lined,boxed,commentsnumbered,vlined,ruled]{algorithm2e}

\DeclareFontFamily{OT1}{pzc}{}
\DeclareFontShape{OT1}{pzc}{m}{it}{<-> s * [1.10] pzcmi7t}{}
\DeclareMathAlphabet{\mathpzc}{OT1}{pzc}{m}{it}

\newtheorem{mydef}{Definition}

\newtheorem{mythm}{Theorem}

\IEEEoverridecommandlockouts                              

\title{\LARGE \bf
Optimal Thermostat Programming and Optimal Electricity Rates for Customers with Demand Charges
}

\author{Reza Kamyar,~\IEEEmembership{Member,~IEEE} and Matthew M. Peet,~\IEEEmembership{Member,~IEEE}
\thanks{*This work was supported by the National Science Foundation under
grant $\#$ CMMI-1301851.}
\thanks{$^{1}$Reza Kamyar and Matthew M. Peet are with the Cybernetic Systems and Controls Lab (CSCL) at the School for Engineering of Matter, Transport and Energy, Arizona State University, Tempe, AZ, 85281, USA,
        {\tt\small rkamyar@asu.edu, mpeet@asu.edu}}%
}

\begin{document}

\maketitle
\thispagestyle{empty}
\pagestyle{empty}

\begin{abstract}
We consider the coupled problems of optimal thermostat programming and optimal pricing of electricity. Our framework consists of a single user and a single provider (a regulated utility). The provider sets prices for the user, who pays for both total energy consumed (\$/kWh, including peak and off-peak rates) and the peak rate of consumption in a month (a demand charge) (\$/kW).
The cost of electricity for the provider is based on a combination of capacity costs (\$/kW) and fuel costs (\$/kWh). In the optimal thermostat programming problem, the user minimizes the amount paid for electricity while staying within a pre-defined temperature range. The user has access to energy storage in the form of thermal capacitance of the interior structure of the building. The provider sets prices designed to minimize the total cost of producing electricity while meeting the needs of the user. To solve the user-problem, we use a variant of dynamic programming. To solve the provider-problem, we use a descent algorithm coupled with our dynamic programming code - yielding optimal on-peak, off-peak and demand prices. We show that thermal storage and optimal thermostat programming can reduce electricity bills using current utility prices from utilities Arizona Public Service (APS) and Salt River Project (SRP). Moreover, we obtain optimal utility prices which lead to significant reductions in the cost of generating electricity and electricity bills.
\end{abstract}


\section{INTRODUCTION}

\begin{figure}[b]%
\begin{center} \vspace{-0.26in}
\subfigure[Percent of growth of electricity demand and its trend-line in the US from 1950 to 2040. Data from~\cite{annual_outlook2014}.]{ \hspace{-0.1in} \includegraphics[scale=0.19]{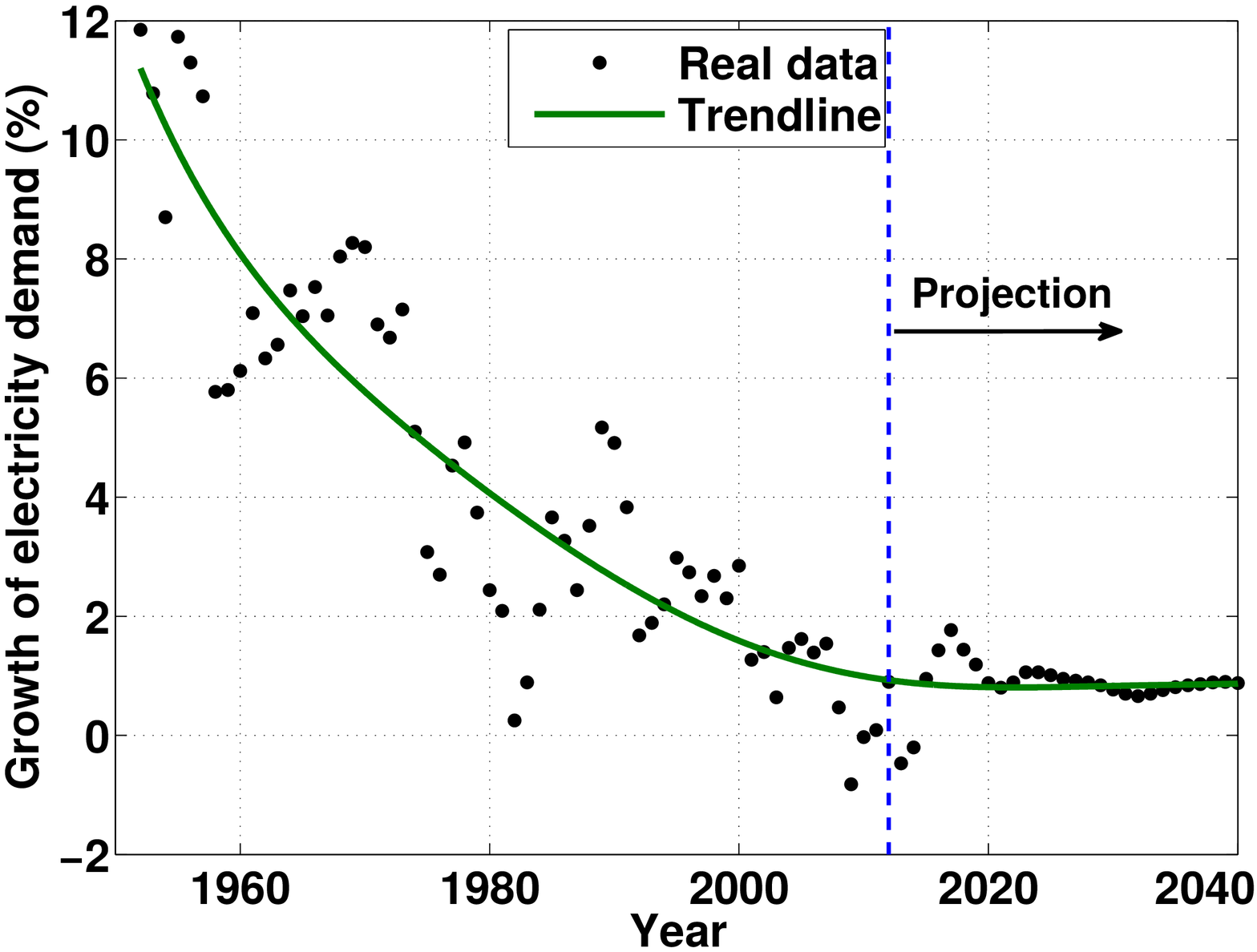}
\label{fig:growth}
} \hspace{0.04in}
\subfigure[Peak-to-average demand of electricity and its trend-line in California and New-England from 1993 to 2012. Data from~\cite{EIA_report}.]{\hspace{-0.15in}\includegraphics[scale=0.19]{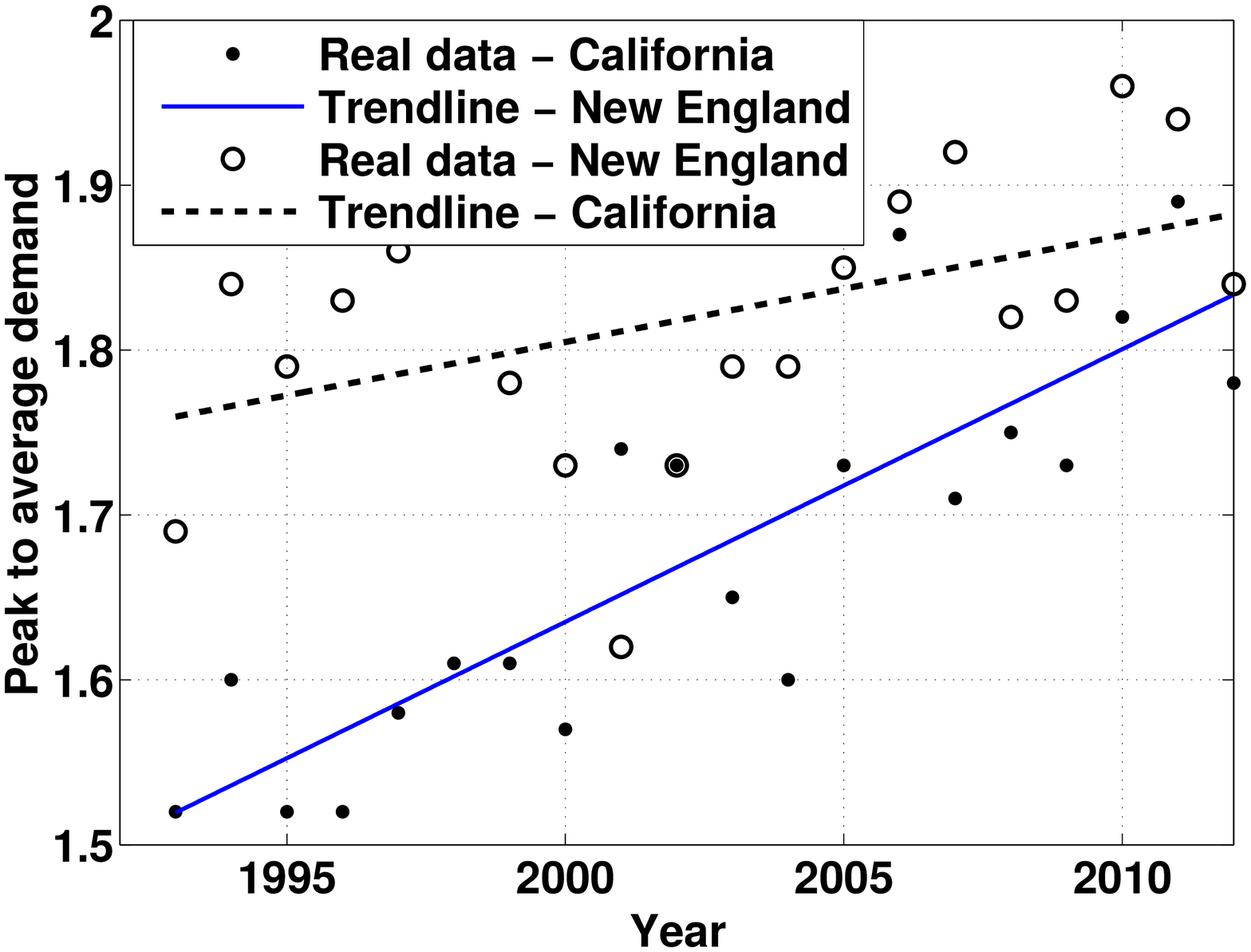}
\label{fig:peak}
}
\end{center} \vspace{-0.15in}
\caption{Demand growth and peak-to-average demand of electricity} \vspace{-0.1in}
\end{figure}

To ensure the reliability of power networks, utility companies must maintain an uninterrupted balance between power generation and demand. In some ways this problem is becoming easier. Partially due to the development of energy-efficient appliances and new materials for insulation, US electricity demand has plateaued~\cite{annual_outlook2014} and is expected to remain flat (less than 1\% growth) for the indefinite future (see Fig.~\ref{fig:growth}). The result is reduced reliance on carbon-producing fossil fuels. However, a new problem has arisen - partially due to increasing use of intermittent renewable energy sources such as distributed solar and wind - in that demand peaks continue to grow. Specifically, as per the US Energy Information Administration (EIA)~\cite{EIA_report}, the ratio of peak demand to average demand has increased dramatically over the last 20 years, setting records of 1.89 in New England in 2012 and 1.96 in California in 2010 (see Fig.~\ref{fig:peak}). Because most utilities are required to maintain generating capacity as determined by peak demand, yet typically only charge customers for total consumption, there is real concern about the viability of existing business models. For example, due to net metering, a typical residential solar customer might have negative consumption during the day and positive consumption during the evening and morning. Such as customer might pay nothing for electricity while contributing substantially to the costs incurred by the utility. In response to this problem, many utilities have sought to halt or reverse growth of the net-metering framework - a process which has met with some limited success.

In this paper, we look at pricing strategies for reducing peak load while retaining the incentives necessary to create a robust distributed renewable sector. Naturally, utilities have been studying this problem for some time and with the widespread adoption of smart-metering (95\% in Arizona), have begun to implement such strategies at scale. Examples of this include on-peak, off-peak and super-peak pricing - rate plans wherein the energy price (\$/kWh) depends on the time of day~\cite{Albdai_review}. By charging more during peak hours, utilities encourage conservation or deferred consumption during hours of peak demand. More aggressive strategies which have emerged recently include voluntary \textit{on-peak demand-limiting} programs wherein customers are rewarded for reducing consumption when requested to do so by the utility~\cite{curtailment}. A yet more aggressive strategy is direct load control~\cite{direct_load_DP,direct_load_DP2} wherein Heating, Ventilating, and Air Conditioning (HVAC) or other appliances are under the direct control of the utilities and can be deferred or deactivated at will. Quite recently, some utilities have introduced demand charges for residential customers. These charges are not based on energy consumption, but rather the maximum \textit{rate of consumption} (\$/kW) over a billing period. While such charges more accurately reflect the cost of generation for the utilities, in practice the effects of such charges on consumption are not well-understood - meaning that the magnitude of the demand charge must be set in an ad-hoc manner (typically proportional to marginal cost of generation).

An alternative approach to reducing peaks in demand is to use energy storage. In this scenario, batteries, pumping and retained heat are used during periods of low demand to create reservoirs of energy which can then be tapped during periods of high demand - thus reducing the need to increase maximum generating capacity. Indeed, the optimal usage of energy storage in a smart-grid environment with dynamic pricing has been recently studied in, for example,~\cite{storage_benefit}. See~\cite{Topcu} for optimal distributed load scheduling in the presence of network capacity constraints. However, to date the high marginal costs of storage infrastructure relative to incentives/marginal cost of additional generating capacity have limited the widespread use of energy storage by consumers/utilities~\cite{battery_usage}. As a cost-free alternative to direct energy storage, it has been demonstrated experimentally~\cite{experiment1,experiment2} and in-silico~\cite{simulation1, simulation2} that the interior structure of buildings and appliances can be exploited as a \textit{passive} thermal energy storage system to reduce the peak-load of the HVAC. A typical strategy - known as \textit{pre-cooling} - is to artificially cool the interior thermal mass (e.g., walls and floor) during periods of low demand. Then, during periods of high demand, heat absorption by these cool interior structures supplements or replaces electricity which would otherwise be consumed by the HVAC. Quantitative assessment of the effect of pre-cooling on demand peak and electricity bills has been evaluated in, e.g.,~\cite{Braun_2006} and {sun2013peak} . It is important to note, however, that ad-hoc strategies such as pre-cooling are only economical when using differential on-peak and off-peak pricing or demand charges.

The goal of this paper is two-fold. First, we consider optimal HVAC usage for a consumer with fixed on-peak, off-peak and demand charges and model passive thermal energy storage using the heat equation. For a given range of acceptable temperatures and using typical data for exterior temperature, we pose the optimal thermostat programming problem as a constrained optimization problem and present a Dynamic Programming (DP) algorithm which is guaranteed to converge to the solution. This yields the temperature set-points which minimize the monthly electricity bill for the consumer. After solving the thermostat programming problem, we use this solution as a model of user behaviour in order to quantify the consumer response to changes in on-peak rates, off-peak rates, and demand charges. We then apply descent methods to this model in order to determine the prices which minimize the cost-of-generation for the utility. In a case study, we show that the optimal prices are NOT necessarily proportional to the marginal costs of generation - meaning that current pricing strategies may be inefficient.

Before presenting our results, we note that models for thermal energy storage do appear in the optimal thermostat programming literature~\cite{curtailment,Braun_2006,guttromson2003residential,henze2004evaluation}. Furthermore, there is an extensive literature on thermostat programming for HVAC systems for on-peak/off-peak pricing~\cite{kelman2011bilinear,lu2005global2, arguello1999nonlinear} as well as \textit{real-time} pricing (prices which are constantly changing)~\cite{constantopoulos1991estia,old_2010,henze2004evaluation,chen2001real} using Model Predictive Control. \cite{Braun_complex_storage} and~\cite{kintner1995optimal} consider optimal thermostat programming with passive thermal energy storage and on-peak/off-peak rates. Perhaps closest to our work, in~\cite{Braun_2006}, the authors use the concept of \textit{deep} and \textit{shallow} mass to create a simplified analogue circuit model of the thermal dynamics of the structure. By using this model and certain assumptions on the gains of the circuit elements,~\cite{Braun_2006} derives an analytical optimal temperature set-point for the demand limiting period which minimizes the demand peak. This scenario would be equivalent to minimizing the demand charge while ignoring on-peak or off-peak rates. Again, referring to~\cite{storage_benefit} and subsequent publications, there has been some excellent work on optimal pricing (albeit without demand charges) for energy storage using batteries in an unregulated electricity market using a social welfare model. This paper differs from existing literature in that it: 1) Considers demand charges (demand charges are far more effective at reducing demand peaks than dynamic pricing) 2) Uses a PDE model for thermal storage (yields a more accurate model of thermal storage) 3) Uses a regulated model for the utility. Although unregulated utility models are popular, the fact is that most US utilities remain regulated.


\section{Problem Statement}
In this section, we first define a model of the thermodynamics which govern heating and cooling of the interior structures of a building. We then use this model to pose the \textit{user-level} (optimal thermostat programming) problem in Sections~\ref{sec:user} as minimization of a monthly electricity bill (with on/peak, off-peak and demand charges) subject to constraints on the interior temperature of the building. Finally, we use this map of on-peak, off-peak and demand prices to consumption to define the \textit{utility-level} problem in Section~\ref{sec:utility} as minimizing the cost of producing electricity.

\subsection{A Model for the Building Thermodynamics}
\label{sec:thermal_mass}

\begin{figure}[ht]
\includegraphics[scale=0.31]{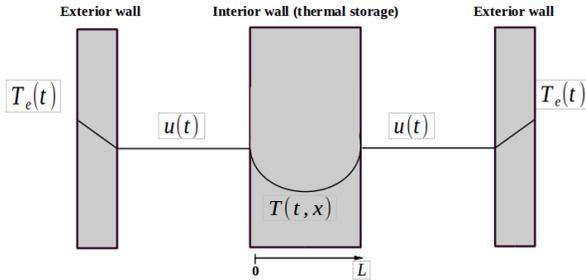} \vspace{-0.2in}
\caption{A schematic view of our thermal mass model}
\label{fig:building}   \vspace{-0.12in}
\end{figure}

To model heat storage in interior walls and floors of a building, we use the one-dimensional unsteady heat conduction equation \vspace{-0.05in}
\begin{equation}
\dfrac{\partial T(t,x)}{\partial t} = \alpha \dfrac{\partial^2 T(t,x)}{\partial x^2}, \vspace{-0.05in}
\label{eq:PDE_conduction}
\end{equation}
where $T: \mathbb{R}^+ \times [0,L_{in}] \rightarrow \mathbb{R}$ represents the temperature distribution in the interior walls/floor with nominal width $L_{in}$ and where $\alpha = \frac{k_{in}}{\rho C_p}$ is the coefficient of thermal diffusivity. Here $k_{in}$ is the coefficient of thermal conductivity, $\rho$ is the density and $C_p$ is the specific heat capacity. The wall is coupled to the interior air temperature using Dirichlet boundary conditions, i.e., $T(t,0) = T(t,L_{in}) = u(t) \;\; \text{for all } t \in \mathbb{R}^+$,
where $u(t)$ represents the interior temperature which we assume can be controlled instantaneously by the thermostat.
 We model the heat loss $q_{loss}$ through the exterior walls by the linear heat sink as \vspace{-0.05in}
\begin{equation}
q_{loss}(t,u(t)) := \dfrac{T_{e}(t)-u(t)}{R_{e}}, \vspace{-0.05in}
\label{eq:qloss}
\end{equation}
where $T_e(t)$ is the outside temperature and $R_{e} = \frac{L_{e}}{k_{e}A_{e}}$ is the thermal resistance of the ext. walls, $L_{e}$ is the nominal width of ext. walls, $k_{e}$ is the coefficient of thermal conductivity and $A_{e}$ is the nominal area of the ext. walls. The heat/energy flux through the surface of the interior wall is modelled as \vspace{-0.05in}
\begin{equation}
q_{in}(T(t,x)) := 2C_{in} \dfrac{\partial T}{\partial x} (t,0), \vspace{-0.05in}
\end{equation}
where $C_{in}=k_{in}A_{in}$ is the thermal capacitance of interior walls and $A_{in}$ is the nominal area of interior walls. By conservation of energy, the power required from the HVAC to maintain the interior air temperature is \vspace{-0.05in}
\begin{equation}
q(t,u(t),T(t,x)) = q_{loss}(u(t),T_e(t)) + q_{in}(T(x,t)).
\label{eq:q}
\end{equation}
See Fig.~\ref{fig:building} for a depiction of the model.

Eqn.~\eqref{eq:PDE_conduction} is a PDE. For optimization purposes, we discretize ~\eqref{eq:PDE_conduction} in space, using $T(t)\in \R^M$ to replace $T(t,x)\in \R$, with $T_i(t)$ denoting $T(t,i\,\Delta x)$, where $\Delta x := \frac{L_{in}}{M+1}$. Then
\begin{equation}
\dot{T}(t) = A \, T(t)+ B \, u(t),
\label{eq:linearized}
\end{equation}
where \vspace{-0.15in}
\[
\begin{small}
A = \dfrac{\alpha}{\Delta x^2}
\begin{pmatrix}
-2 & 1 & 0 & 0  \\
1 & \ddots & \ddots & 0  \\
0 & \ddots & \ddots & 1  \\
0 & 0 & 1 & -2
\end{pmatrix}  , \quad
B= \dfrac{\alpha}{\Delta x^2}
\begin{pmatrix}
1 \\
0 \\
\vdots \\
0 \\
1
\end{pmatrix}\in \mathbb{R}^M.
\end{small}
\]
By discretizing in time, using $ \dot{T}(t) \approx \left( T(t+\Delta t) - T(t) \right)/ \Delta t $ we can rewrite Equation~\eqref{eq:linearized} as a difference equation. 
\begin{small}
\begin{align}
T^{k+1}  =\hspace{-0.03in}
\begin{bmatrix}
T^{k+1}_{1} \\
\vdots \\
T^{k+1}_{M}
\end{bmatrix} \hspace{-0.03in} = \hspace{-0.03in}
 f(T^{k},u_{k})
& \hspace{-0.03in} =  \hspace{-0.03in} \begin{bmatrix}
f_1(T^{k},u_{k}) \\
\vdots \\
f_M(T^{k},u_{k})
\end{bmatrix} \hspace{-0.02in}
 =  (I + A \, \Delta \, t )T^{k} + B \, \Delta t \, u_{k}
\label{eq:discrete_dyn}
\end{align}
\end{small}\hspace*{-0.05in}
for $k=0, \cdots, N_f-1$, where  now we have $T^k=T(k\, \Delta t)$ and $u_k=u(k \, \Delta t)$. \vspace{-0.1in}

\subsection{User-Level Problem: Optimal Thermostat Programming}
\label{sec:user}
In this section, we define the problem of optimal thermostat programming. We first divide each day into three periods: off-peak hours from 12 AM to $t_{\text{on}}$ with electricity price $p_{\text{off}} \, (\$/kWh)$; on-peak hours beginning at $t_{\text{on}}$ and ending at $t_{\text{off}} > t_{\text{on}}$ with electricity price $p_{\text{on}}  \,(\$/kWh)$; and off-peak hours from $t_{\text{off}}$ to 12 AM with electricity price $p_{\text{off}} \, (\$/kWh)$. In addition to the on-peak and off-peak charges, we consider a monthly charge which is proportional to the maximum rate of consumption during the peak hours. The proportionality constant is called the \textit{demand price} $p_d \, (\$/kW)$. Given prices $p_{\text{on}}, p_{\text{off}}$ and $p_d$, the total cost of consumption (daily electricity bill) is divided as
\begin{equation}
J_t(\mathbf{u} ,T_1, p_{\text{on}},p_{\text{off}},p_d) = J_e(\mathbf{u} ,T_1,p_{\text{on}},p_{\text{off}}) + J_d(\mathbf{u} ,T_1, p_d),
\label{eq:Jt}
\end{equation}
where $J_e$ is the energy cost, $J_d$ is the demand cost and
\[
\mathbf{u} := [u_0, \cdots, u_{N_f-1}] \in \mathbb{R}^{N_f}.
\]
The energy cost is
\begin{align}
& \hspace{-0.1in} J_e(\mathbf{u}, T_1,p_{\text{on}},p_{\text{off}}) \nonumber \\
& \quad = \left( p_{\text{off}} \sum_{k \in S_{\text{off}}} g(k,u_{k},T_1^k)  + p_{\text{on}}  \sum_{k \in S_{\text{on}}} g(k,u_{k},T_1^k) \right) \Delta t,
\label{eq:Je}
\end{align}
where $k \in S_{\text{on}}$ if $k \Delta t \in [t_{on},t_{off}]$ and $k\in S_{\text{off}}$ otherwise. That is, $S_{\text{on}}$ and $S_{\text{off}}$ correspond to the set of on-peak and off-peak sampling times, respectively. The function $g$ is a discretized version of $q$ (Eqn.~\eqref{eq:q}):
\begin{align}
\hspace{-0.04in} g(k,u_k, T_1^k) := \dfrac{T_e(k \, \Delta t) - u(k \, \Delta t)}{R_e}& +  2 \, C_{in} \dfrac{T(k \, \Delta t,\Delta x) - u(k \, \Delta t)}{\Delta x} \nonumber \\
 & \hspace{-0.25in} =
 \dfrac{T^k_{e} - u_k}{R_e} + 2 \, C_{in} \dfrac{T_1^k - u_k}{\Delta x}.
 \label{eq:g}
\end{align}
This is the power consumed by the HVAC, where $T_e^k$ denotes the external temperature at time-step $k$. If demand charges are calculated monthly, the demand cost $J_d$ for a single day~is
\begin{equation}
J_d(\mathbf{u}, T_1, p_d) := \dfrac{p_d}{30} \max_{k \in S_{\text{on}}} g(k,u_{k},T_1^k).
\label{eq:Jd}
\end{equation}

We now define the optimal thermostat programming problem at the user-level as minimization of the total cost of consumption $J_t$ as defined in~\eqref{eq:Jt}, subject to the building thermodynamics (Eqn.~\eqref{eq:discrete_dyn}) and interior temperature constraints ($T(t) \in [T_{\min},T_{\max}]$). \vspace{-0.05in}
\begin{align}
& J^*(p_{\text{on}},p_{\text{off}},p_d) = \min_{u_{k}, \gamma \in \mathbb{R}, T^k \in \mathbb{R}^M} J_{e}(\mathbf{u}, T_1,p_{\text{on}},p_{\text{off}}) + \frac{p_d}{30} \, \gamma   \nonumber \\
& \text{subject to} \nonumber  \;\, g(k, u_{k}, T_1^k) \leq \gamma \hspace{0.92in} \text{ for } k \in S_{\text{on}} \nonumber \\
& \hspace{0.58in} T^{k+1} = f(T^k,u_k) \hspace{0.83in} \text{ for } k \in S_{\text{on}} \cup S_{\text{off}} \nonumber \\
& \hspace{0.58in}  T_{\min} \leq u_{k} \leq T_{\max} \hspace{0.85in}   \text{ for } k \in S_{\text{on}} \cup S_{\text{off}}  \nonumber \\
&  \hspace{0.58in} T^{0} = [T_{\text{init}}(\Delta x), \cdots, T_{\text{init}}(M \, \Delta x)]^T,
\label{eq:user_discrete}
\end{align}
where $T_{\min}$ and $T_{\max}$ are the acceptable bounds on the interior temperature. Note that this optimization problem depends implicitly on exterior temperature through the time-varying function $g$.

\subsection{Utility-Level Optimization Problem}
\label{sec:utility}

We define the utility-level optimization problem as minimization of the cost of generating electricity such that generation is equal to consumption, and revenue is equal to cost of generation. Let $s(t)$ be the amount of electricity produced as a function of time and $s_k=s(k\Delta t)$. 
First, we consider a linear model of the production cost (adopted from Arizona Public Utility SRP) as
\[
a \hspace{-0.1in} \sum_{k \in S_{\text{on}} \cup S_{\text{off}}} \hspace{-0.1in} s_k + b\max_{k\in S_{\text{on}}} s_k ,
\]
where $a  \, (\$/kWh)$ is the marginal cost of producing the next $kWh$ of energy and $b \, (\$/kW)$ is the marginal cost of installing the next $kW$ of production capacity. Values of the coefficients $a$ and $b$ for SRP can be found in e.g.,~\cite{SRP}. Now define $\mathbf{u}^*(p_{\text{on}},p_{\text{off}},p_d)$ and $T^*(p_{\text{on}},p_{\text{off}},p_d)$ to be minimizing arguments to the user-level problem defined in~\eqref{eq:user_discrete}. Then the constraint that production equals consumption implies $s_k = g(k,u^*_{k}(p_{\text{on}},p_{\text{off}},p_d),T_1^{*k}(p_{\text{on}},p_{\text{off}},p_d))$.
We now define the utility-level optimization problem as minimization of the cost of electricity production subject to equality of production and consumption.
\begin{align}
& c^* :=  \min_{ p_{\text{on}},p_{\text{off}},p_d \in \R} \quad  a \hspace{-0.1in} \sum_{k \in S_{\text{on}} \cup S_{\text{off}}} \hspace{-0.1in} s_k + b\max_{k\in S_{\text{on}}} s_k  \label{eq:utility_problem} \\
& \text{subject to } \; s_k = g(k,u^*_{k}(p_{\text{on}},p_{\text{off}},p_d),T_1^{*k}(p_{\text{on}},p_{\text{off}},p_d)),  \nonumber\\
& \qquad \qquad \;\;\, a \hspace{-0.05in} \sum_{k \in S_{\text{on}} \cup S_{\text{off}}} \hspace{-0.05in} s_k + b \, \max_{k\in S_{\text{on}}} s_k \, = \, J_t(\mathbf{u}^*(p_{\text{on}},p_{\text{off}},p_d),   \nonumber \\ 
& \hspace{1.44in} T^*_1(p_{\text{on}},p_{\text{off}},p_d), p_{\text{on}},p_{\text{off}},p_d),
 \nonumber
\end{align}
where the last two lines constrain that costs equal revenue (recall $J_t$ is revenue from the users as defined in~\eqref{eq:Jt}).


\section{Solving User- and Utility-level Problems}

First, we solve the optimal thermostat programming problem using a variant of dynamic programming. This yields consumption as a function of prices $p_{\text{on}},p_{\text{off}},p_d$. Next, we embed this implicit function in a descent algorithm in order to find prices which minimize the Utility-level optimization problem as formulated in~\eqref{eq:utility_problem}. We start by defining a \textit{cost-to-go} function, $V_k$. Given $\gamma \in \mathbb{R}^+$, at the final time $N_f \, \Delta t=24$, we have
\begin{equation}
V_{N_f}(x) := \dfrac{p_d}{30} \gamma.
\label{eq:V_Nf}
\end{equation}
Here for simplicity, we use $x=T \in \R^M$ to represent the discretized temperature distribution in the wall. Define prices $p_{j}=p_{\text{off}}$ if $j\in S_{\text{off}}$ and $p_{j}=p_{\text{off}}$ otherwise. Then, we construct the cost-to-go function inductively as
\begin{align}
\hspace{-0.1in} V_{j-1} (x) := \hspace{-0.1in}  \min_{u \in W_{\gamma,j-1}(x)}  \hspace{-0.05in}  \left( p_{j-1} \, g(j-1, u,x_1) \Delta t + V_j(f(x,u)) \right),
\label{eq:V0_on}
\end{align}
where $W_{\gamma,j}(x)$ is the set of allowable inputs at time $j$ and state $x$: \vspace{-0.1in}
\begin{align*}
&W_{\gamma,j}(x) := \\
&\begin{cases} \{u \in \mathbb{R}: T_{\min} \leq u \leq T_{\max}, g(j-1,u,x_1) \leq \gamma\} & j \in S_{\text{on}}\\
\{u \in \mathbb{R}: T_{\min} \leq u \leq T_{\max} \}& j \in S_{\text{off}}
\end{cases}.
\end{align*}
Now we present the main result.
\begin{mythm} Given $\gamma \in \mathbb{R}^+$, suppose that $V_i$ satisfies~\eqref{eq:V_Nf} and~\eqref{eq:V0_on}. Then $V_0(T^0)=J^*$, where
\begin{align}
& J^*(p_{\text{on}},p_{\text{off}},p_d) = \min_{u_{k}, T^k \in \mathbb{R}^M} J_{e}(\mathbf{u}, T_1,p_{\text{on}},p_{\text{off}}) + \frac{p_d}{30} \, \gamma   \nonumber \\
& \text{subject to} \nonumber  \;\, g(k, u_{k}, T_1^k) \leq \gamma \hspace{0.92in} \text{ for } k \in S_{\text{on}} \nonumber \\
& \hspace{0.58in} T^{k+1} = f(T^k,u_k) \hspace{0.83in} \text{ for } k \in S_{\text{on}} \cup S_{\text{off}} \nonumber \\
& \hspace{0.58in}  T_{\min} \leq u_{k} \leq T_{\max} \hspace{0.85in}   \text{ for } k \in S_{\text{on}} \cup S_{\text{off}}  \nonumber \\
&  \hspace{0.58in} T^{0} = [T_{\text{init}}(\Delta x), \cdots, T_{\text{init}}(M \, \Delta x)]^T.
\label{eq:user_discrete2}
\end{align}
\label{thm:DP} \vspace{-0.1in}
\end{mythm}

To prove Theorem~\ref{thm:DP}, we require the following definitions.  
\begin{mydef}
Given $p_{\text{off}}, p_{\text{on}}, p_d, \gamma \in \mathbb{R}^+$, $N_f \in \mathbb{N}^+$, and $t_{\text{off}},t_{\text{on}},\Delta t \in \mathbb{R}^+$ such that $\frac{t_{\text{on}}}{\Delta t}, \frac{t_{\text{off}}}{\Delta t} \in \mathbb{N}$, define the cost-to-go functions
\[
Q_j: \mathbb{R}^{N_f-j} \times \mathbb{R}^{N_f-j+1} \rightarrow \mathbb{R} \; \text{ for } j=0, \cdots, N_f \text{ as} \vspace{-0.1in}
\]
\vspace{-0.1in}
\begin{small}
\begin{align}
& \hspace{-3.3in} Q_j(x, y, p_{\text{on}}, p_{\text{off}}) := \nonumber \\
\begin{cases}
    \!\begin{aligned}
       & \hspace{-0.03in} \left( p_{\text{off}} \hspace{-0.2in} \sum\limits_{\substack{k \in S_{\text{off}} \\ k \notin \{ 0, \cdots, j-1\}}} \hspace{-0.2in} g(k,x_{k},y_{k}) + p_{\text{on}} \hspace{-0.05in} \sum\limits_{k \in S_{\text{on}}} \hspace{-0.05in} g(k,x_{k},y_{k}) \hspace{-0.02in} \right) \Delta t & \hspace{-0.05in}
    \end{aligned}   & \hspace{-0.1in} \text{ if } 0 \leq j < N_{\text{on}}  \\
\hspace{-0.03in} \left(  p_{\text{on}} \hspace{-0.2in} \sum\limits_{\substack{k \in S_{\text{on}} \\ k \notin \{ N_{\text{on}},\cdots,j-1\} }} \hspace{-0.28in} g(k,x_{k},y_{k}) + p_{\text{off}} \hspace{-0.27in} \sum\limits_{\substack{k \in S_{\text{off}} \\ k \notin \{0,\cdots,N_{\text{on}}-1 \}}} \hspace{-0.28in} g(k,x_{k},y_{k}) \hspace{-0.03in} \right) \Delta t  & \hspace{-0.1in} \text{ if } N_{\text{on}} \leq j < N_{\text{off}} \\
\;\; p_{\text{off}} \hspace{-0.2in}  \sum\limits_{\substack{ k \in \{j,\cdots, N_f-1 \}}} \hspace{-0.2in} g(k,x_{k},y_{k})  \Delta t   & \hspace{-0.1in} \text{ if } N_{\text{off}} \leq j < N_f \\
\dfrac{p_d}{30} \gamma  & \hspace{-0.1in} \text{ if } j = N_f,
\end{cases}
\label{eq:cost-to-go} 
\end{align}
\end{small}
\vskip-0.35cm
\belowdisplayskip=0pt
\hspace{-0.18in} where $g$ is defined as in~\eqref{eq:g}, and $N_{\text{on}} := \frac{t_{\text{on}}}{\Delta t}$ and $N_{\text{off}} := \frac{t_{\text{off}}}{\Delta t}$ are the time-steps corresponding to start and end of the on-peak hours.
\label{def:cost-to-go}
\end{mydef}

Note that from~\eqref{eq:Je}, it is clear that $Q_0 = J_e + \dfrac{p_d}{30} \gamma$. \vspace{0.1in}

\begin{mydef}
\hspace{-0.05in}Given $\gamma, T_{\min}, T_{\max} \in \mathbb{R}$ and $N_f,M \in \mathbb{N}^+$, define the set
\begin{align}
U_j(x) &:= \{ (u_j,  \cdots,  u_{N_f-1}) \in \mathbb{R}^{N_f-j} : \nonumber \\
& \hspace{-0.08in}  g(k,u_{k}, T_1^{k}) \leq \gamma \; \text{ for all }   k \in S_{\text{on}}, \nonumber \\
& \hspace{-0.08in}  T^j = x \text{ and } T^{k+1} = f(T^{k},u_{k})  \text{ for all } k \in \{ j, \cdots, N_f-1\}, \nonumber \\
& \hspace{-0.08in}   T_{\min} \leq u_{k} \leq T_{\max}  \text{ for all }   k \in S_{\text{on}} \cup S_{\text{off}} \}
\label{eq:Uj}
\end{align}
for any $x \in \mathbb{R}^M$ and for every $j \in \{ 0, \cdots, N_f-1 \}$, where $f$ and $g$ are defined as in~\eqref{eq:discrete_dyn} and~\eqref{eq:g}.  \vspace{0.1in}
\label{def:Uj}
\end{mydef} 

\begin{mydef}
Given $N_f, M \in \mathbb{N}^+$, $j \in \{ 0, \cdots, N_f-1 \}$, let
\[
\overline{\mu}_{j} :=[ \mu_j, \cdots, \mu_{N_f-1} ]
\]
where $ \mu_k : \mathbb{R}^{M} \rightarrow \mathbb{R} \text{ for } k=j, \cdots, N_f-1.
$
Consider $U_j$ as defined in~\eqref{eq:Uj} and $f$ as defined in~\eqref{eq:discrete_dyn}. If
\[
\overline{\mu}_j(w) := [\mu_j(w), \mu_{j+1}(T^{j+1}) \cdots, \mu_{N_f-1}(T^{N_f-1})] \in U_j(T^j)
\]
for any $w \in \mathbb{R}^M$, where 
\[
T^{k+1} = f(T^k,\mu_k(T^k)), T^j=w \;\; \text{for } k=j ,\cdots, N_f-2,
\]
then we call $\overline{\mu}_{j}$  an \textit{admissible control law} for the system
\[T^{k+1} = f(T^k,\mu_k(T^k)), \; k=j ,\cdots, N_f-1
\]
for any $w \in \mathbb{R}^M$.
\label{def:maps}
\end{mydef}

We now present a proof for Theorem~\ref{thm:DP}.

\begin{proof}
Since the cost-to-go function $Q_0 = J_e + \dfrac{p_d}{30} \gamma$, if we show that \vspace{-0.08in}
\begin{equation}
\min_{ \overline{\mu}_j(T^j) \in U_j(T^j)} \; Q_j(\overline{\mu}_{j}(T^j),T_1, p_{\text{on}}, p_{\text{off}}) = V_{j}(T^j) \vspace{-0.05in}
\label{eq:proof1}
\end{equation}
for $j=0, \cdots, N_f$ and for any $T^j \in \mathbb{R}^{M}$, where
\[
T_1:=[T^j,f(T^j,\mu_j(T^j)), \cdots, f(T^{N_f-1},\mu_{N_f-1}(T^{N_f-1}))],
\] 
then it will follow that $J^* = V_0(T^0)$. For brevity, we denote $\overline{\mu}_j(T^j)$ by $\overline{\mu}_j$, $U_j(T^j)$ by $U_j$ and we drop the last two arguments of $Q_j$.
To show~\eqref{eq:proof1}, we use induction as follows. 

\noindent \textit{\underline{Basis step}}: If $j=N_f$, then from~\eqref{eq:V_Nf} and~\eqref{eq:cost-to-go} we have $V_{N_f}(T^{N_f}) = \frac{p_d}{30} \gamma$.\vspace*{0.05in}

\noindent \textit{\underline{Induction hypothesis}}: Suppose \vspace{-0.065in}
\[
\min_{\overline{\mu}_k \in U_k} \; Q_k(\overline{\mu}_{k}, T_1) = V_{k}(T^k) \vspace{-0.065in}
\]
for some $k \in \{0, \cdots, N_f \}$ and for any $T^k \in \mathbb{R}^{M}$. Then, we need to prove that \vspace{-0.075in}
\begin{equation}
\min_{ \overline{\mu}_{k-1} \in U_{k-1}} \; Q_{k-1}(\overline{\mu}_{k-1},T_1) = V_{k-1}(T^{k-1}) \vspace{-0.05in}
\label{eq:proof4}
\end{equation}
for any $T^k \in \mathbb{R}^{M}$. Here, we only prove~\eqref{eq:proof4} for the case which $N_{\text{off}} < k \leq N_f-1$.
The proofs for the cases $0 \leq k \leq N_{\text{on}}$ and $N_{\text{on}} < k \leq N_{\text{off}}$ follow the same exact logic.

Assume that $N_{\text{off}} < k \leq N_f-1$. Then, from Definition~\ref{def:cost-to-go} 
\begin{align}
& \hspace{-0.1in}\min_{ \overline{\mu}_{k-1} \in U_{k-1}} \;  Q_{k-1}(\overline{\mu}_{k-1},T_1) \nonumber  \\
=& \hspace{-0.1in}  \min_{\mu_{k-1}, \cdots, \mu_{N_f-1} \in R} \; p_{\text{off}} \left( \sum\limits_{j=k-1}^{N_f-2}  g(j,\mu_j,T_1^{j}) \right) \Delta t \nonumber  \\
=& \hspace{-0.1in} \min_{\mu_{k-1}, \cdots, \mu_{N_f-1} \in R}  p_{\text{off}} \left( g(k-1, \mu_{k-1},T_1^{k-1}) \right.
 \left. \hspace{-0.05in} + \hspace{-0.05in} \sum\limits_{j=k}^{N_{f}-2}  g(j, \mu_j ,T_1^{j}) \hspace{-0.035in} \right) \Delta t,
\label{eq:proof5}
\end{align}
\vskip-0.08in
\hspace{-0.18in} where $R:= \{ x \in \mathbb{R}: T_{\min} \leq x  \leq T_{\max} \}$.
From the principle of optimality~\cite{bellman} it follows that \vspace{-0.07in}
\begin{align}
&\hspace{-0.15in}\min_{\mu_{k-1}, \cdots, \mu_{N_f-1} \in R} \hspace{-0.01in}   p_{\text{off}} \left( \hspace{-0.03in} g( k-1,\mu_{k-1},T_1^{k-1} \hspace{-0.01in}   ) + \hspace{-0.05in} \sum\limits_{j=k}^{N_f-1} \hspace{-0.04in}  g(j, \mu_j ,T_1^{j}) \hspace{-0.035in} \right) \Delta t \nonumber \\
=&\min_{\mu_{k-1} \in R} \; \left( p_{\text{off}} \, g(k-1,\mu_{k-1},T_1^{k-1})  \Delta t \right. \nonumber  \\
 & \left. \hspace{0.6in} + \min_{\mu_{k}, \cdots, \mu_{N_f-1} \in R} \;  p_{\text{off}} \sum\limits_{j=k}^{N_f-1}  g( j,\mu_{j},T_1^{j}) \right) \Delta t,
 \label{eq:proof6}
\end{align}
By combining~\eqref{eq:proof5} and~\eqref{eq:proof6} we have \vspace{-0.07in}
\begin{align}
& \min_{\overline{\mu}_{k-1} \in U_{k-1}} \; Q_{k-1}(\overline{\mu}_{k-1}, T_1) \nonumber \\
&= \min_{\mu_{k-1} \in R} \; \left( p_{\text{off}} \, g(k-1,\mu_{k-1}),T_1^{k-1})  \Delta t \right. \nonumber  \\
 &   \hspace{0.7in} \left. + \min_{\mu_{k}, \cdots, \mu_{N_f-1} \in R} \;  p_{\text{off}} \sum\limits_{j=k}^{N_f-1}  g(j, \mu_{j},T_1^{j}) \right) \Delta t.
 \label{eq:proof7}
\end{align}
\vskip-0.1in
\hspace{-0.18in}
From Definition~\ref{def:cost-to-go}, we can write \vspace{-0.05in}
\begin{align}
\min_{\mu_{k}, \cdots, \mu_{N_f-1}} \hspace{-0.04in}  p_{\text{off}} \, \Delta t \sum\limits_{j=k}^{N_f-1}  g(j,\mu_{j},T_1^{j})
= \hspace{-0.07in} \min_{\overline{\mu}_{k} \in U_k}   Q_{k}(\overline{\mu}_{k},T_1).
\label{eq:proof2}
\end{align}
Then, by combining~\eqref{eq:proof7} and~\eqref{eq:proof2} and using the induction hypothesis it follows that
\begin{align*}
& \min_{\overline{\mu}_{k-1} \in U_{k-1}}  \; Q_{k-1}(\overline{\mu}_{k-1},T_1) \\
&  = \min_{\mu_{k-1} \in R} \left( p_{\text{off}} \, g(k-1,\mu_{k-1},T_1^{k-1}) \Delta t  + \hspace{-0.05in}  \min_{\overline{\mu}_{k} \in U_k} Q_{k}(\overline{\mu}_{k},T_1) \right) \\
& \hspace{0.7in} = \min_{\mu_{k-1} \in R} \; \left( p_{\text{off}} \, g(k-1,\mu_{k-1},T_1^{k-1}) \Delta t  + V_{k}(T^k) \right)
\end{align*}
for any $T^k \in \mathbb{R}^{M}$.
By substituting for $T^k$ from~\eqref{eq:discrete_dyn} and using the definition of $V$ in~\eqref{eq:V0_on} we have
\begin{align*}
\min_{\overline{\mu}_{k-1} \in U_{k-1}} \hspace{-0.05in} Q_{k-1} ( & \overline{\mu}_{k-1}, T_1)   = 
\min_{\mu_{k-1} \in R} \left( p_{\text{off}} \, g(k-1, \mu_{k-1},T_1^{k-1}) \Delta t \right. \\
&  \quad\left. + V_{k}(f(T^{k-1},\mu_{k-1}(T^{k-1}))) \right) = V_{k-1}(T^{k-1})
\end{align*}
for any $T^{k-1} \in \mathbb{R}^M$.
By using the same logic it can be shown that
$
\min_{\overline{\mu}_{k-1} \in U_{k-1}} Q_{k-1}(\overline{\mu}_{k-1},T_1) = V_{k-1}(T^{k-1})
$
for any $k \in \{0, \cdots, N_{\text{off}}-1 \}$ and for any $T^{k-1} \in \mathbb{R}^{M}$.
Therefore, by induction,~\eqref{eq:proof1} is true. Thus, $J^* = V_0(T^0)$.
\end{proof}

Using Theorem~\ref{thm:DP}, we propose Algorithm~\ref{algorithm} to find solutions to the user-level problem~\eqref{eq:user_discrete} and the utility-level problem~\eqref{eq:utility_problem}.  \vskip-0.12in

\begin{algorithm}
\begin{footnotesize}
\textbf{Inputs:}\\
External temperature $T_e$, start and end of on-peak hours $N_{\text{on}}, N_{\text{off}}$, thermal resistance $R_e$, thermal capacitance $C_{in}$, initial temperature $T_{\text{init}}$ of walls, step-sizes $\Delta t$ and $\Delta x$, minimum and maximum interior temperatures $T_{\min}, T_{\max}$, marginal costs $a$ and $b$, step-sizes $\Delta p_d$ and $\Delta p_{\text{on}}$ on electricity prices, initial prices $p_{d_0}$ and $p_{\text{on}_0}$ ($p_{d_0} + p_{\text{on}_0} < 1$), maximum number of bisection iterations $b_{\max}$, lower bound $\gamma_l$ and upper bound $\gamma_u$ for bisection search, stopping threshold $\epsilon$.

\vspace{0.05in}

\textbf{Initialization:} \\
Set $p_d = p_{d_0}, \, p_{\text{on}} = p_{\text{on}_0}, \, p_{\text{off}} = 1 - p_d -p_{\text{on}}$. 
	Set $k=0$. \\
	\While{ $k \leq b_{\max}$}{	
		Set $\gamma = \frac{\gamma_u+\gamma_l}{2}$. 	\\	
  		\eIf{ $V_0$ in~\eqref{eq:V0_on} exists }{		
  			Calculate $u_0^*, \cdots, u_{N_f-1}^*$ as the minimizers of the RHS of~\eqref{eq:V0_on}. 	
  			Set $\gamma_u = \gamma$.				
 		}
 		{
	 	 	Set $\gamma_l = \gamma$. \\
 		}					
		Set $k = k+1$.
	}
Calculate $F_{old} = a\, G+b \, g_{\max}$ as defined in~\eqref{eq:utility_problem}. Set $F_{new} = F_{old} + 2 \, \epsilon$. \vspace{0.05in}\\

\textbf{Main loop:}\\

\While{ $F_{\text{new}} - F_{\text{old}} > \epsilon $}{

Set $F_{old} = F_{new}$.

\For{ $s_d  \in \{ -\Delta p_d, \Delta p_d \}$}{
	\For{ $s_{\text{on}} \in \{ -\Delta p_{\text{on}}, \Delta p_{\text{on}} \}$}{
	
	Set $p_{d} = p_d + s_d, p_{\text{on}} = p_{\text{on}} + s_{\text{on}}, p_{\text{off}} = 1 - p_d -p_{\text{on}}$.\\
	Set $k=0$. \\
	\While{ $k \leq b_{\max}$}{	
		Set $\gamma = \frac{\gamma_u+\gamma_l}{2}$. 	\\	
  		\eIf{ $V_0$ in~\eqref{eq:V0_on} exists }{		
  			Calculate $u_0^*, \cdots, u_{N_f-1}^*$ as the minimizers of the RHS of~\eqref{eq:V0_on}. 	
  			Set $\gamma_u = \gamma$.				
 		}
 		{
	 	 	Set $\gamma_l = \gamma$. \\
 		}					
		Set $k = k+1$.
	}
	Calculate $ \text{cost} = a\, G+b \, g_{\max}$ as defined in~\eqref{eq:utility_problem}.\\
	\If{ $\text{cost} \leq F_{\text{new}}$}{		
  			Set $F_{new} = cost$. 	
  			Set 	$\mathbf{u}^* = [u_0^*, \cdots, u_{N_f-1}^*]$.  \\
  			Set $J^*_t=J_t(\mathbf{u}^*, T_1, p_{\text{on}}, p_{\text{off}}, p_d)$, $J_t$ defined in~\eqref{eq:Jt}-\eqref{eq:Jd}.
  			Set $p^* = \frac{\text{cost}}{J^*_t} \cdot [p_{\text{on}}, p_{\text{off}}, p_d]$.
 		}
	}
}
}

\textbf{Outputs:}
 Optimal prices $p^*$ and optimal interior temperature $\mathbf{u}^*$.
\end{footnotesize}
\caption{A descent algorithm for computing optimal electricity prices}
\label{algorithm}
\end{algorithm}
\vspace{-0.15in}

\section{Numerical Examples and Analysis}
\label{sec:numerical}
\begin{figure}[htbp]
\vspace{-0.2in}
\centering
\includegraphics[scale=0.3]{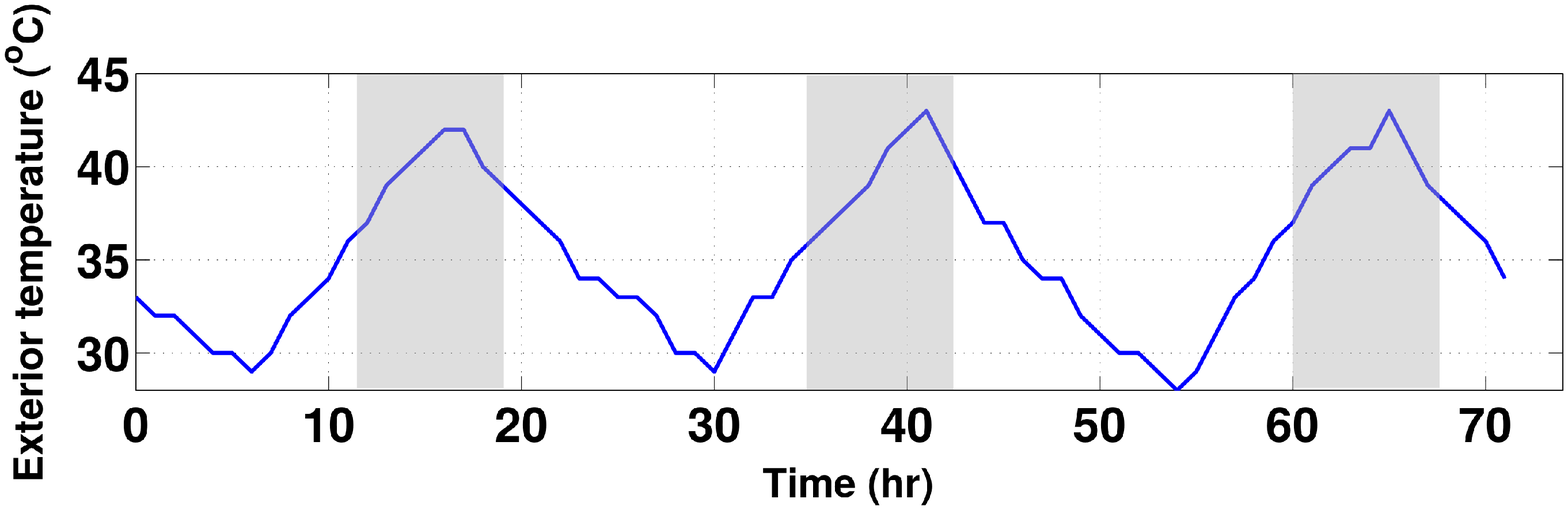} \vspace{-0.15in}
\caption{External temperature of three typical summer days in Phoenix, Arizona. Shaded areas correspond to on-peak hours.}
\label{fig:Te} 
\end{figure}

In this section, we demonstrate convergence of our algorithm for optimal thermostat programming using electricity prices from APS and temperature data from Phoenix, AZ. In addition, we study the problem of optimal electricity pricing using the marginal cost data from SRP. We ran all the numerical simulations for three consecutive days with time-step $\Delta t = 1 \; hr$, space-step $\Delta x = 0.1 \; m$ and with building's parameters listed in Table~\ref{tab:params}.
We used an external temperature profile measured for three typical summer days in Phoenix, Arizona (see Fig.~\ref{fig:Te}). For each day, the on-peak period starts at $N_{\text{on}} = 12 $ PM and ends at $N_{\text{off}} = 7$ PM. In all scenarios, we used $T_{\min} = 22^{\circ}C$ and $T_{\max} = 28^{\circ}C$.

\begin{table}[htbp]
\caption{Building's parameters as defined in Section~\ref{sec:thermal_mass}}
\vspace{-0.1in}
\centering
\begin{tabular}{|c|c|c|c|}
\hline
$L_{in} (m)$ & $\alpha (m^2/s)$ & $R_e (K/W)$ & $C_{in} (Wm/K)$ \\
\hline
0.4 & $8.3 \times 10^{-7}$ & 0.0015 & 45  \\
\hline
\end{tabular}
\label{tab:params}
\end{table}
\vspace{-0.2in}

\begin{table}[htbp]
\renewcommand{\tabcolsep}{2.5pt}
\caption{on-peak, off-peak and demand prices, Arizona utility APS~\cite{APS}}
\vspace{-0.1in}
\centering
\begin{tabular}{|c|c|c|c|}
\hline
 & on-peak ($\$$ per $kWh$) & off-peak ($\$$ per $kWh$) & demand ($\$$ per $kW$)  \\
\hline
APS & 0.089  & 0.044  & 13.50  \\
\hline
\end{tabular}
\label{tab:rates} \vspace{-0.15in}
\end{table}

\subsection{Scenario 1: The Effect of Electricity Prices on Peak Demand
and Production Costs}
\begin{figure}[htbp]
\includegraphics[scale=0.345]{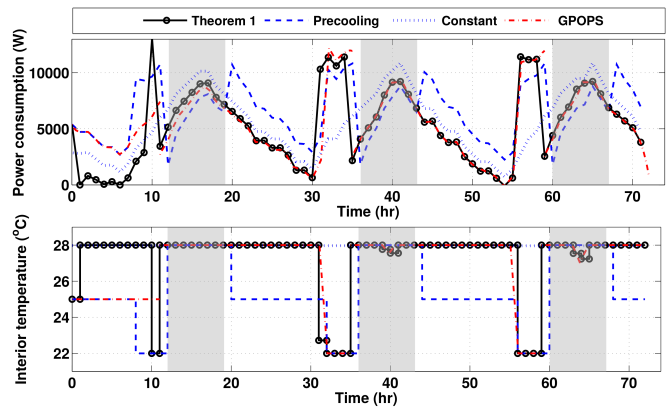}
\caption{Utility power consumption and temperature settings for various programming strategies using APS's rates.}
\label{fig:APS}
\end{figure}

In this scenario, we first consider the optimal thermostat programming problem (See \eqref{eq:user_discrete}) using the electricity prices $p_{\text{off}}, p_{\text{on}}$ and $p_d$ as determined by APS~\cite{APS} (See Table~\ref{tab:rates}). The results of the dynamic programming algorithm are given in Table~\ref{tab:APS} as total price paid by the user (we also include the peak demand). For comparison, we have run the same optimal control problem using the general-purpose optimization solver GPOPS~\cite{gpops}. Moreover, we have compared our result with a typical precooling strategy and a naive strategy of setting the temperature to $T_{\max}$ (constant). As can be seen, our algorithm outperforms the heuristic approaches. The power consumption and the temperature setting as a function of time for each strategy can be found in Fig.~\ref{fig:APS}. For convenience, the on-peak and off-peak intervals are indicated on the figure. As can be seen, for APS prices and our building's parameters, the optimal strategy does not reduce the peak demand with respect to the precooling strategy.

\begin{figure}[htbp]
\includegraphics[scale=0.33]{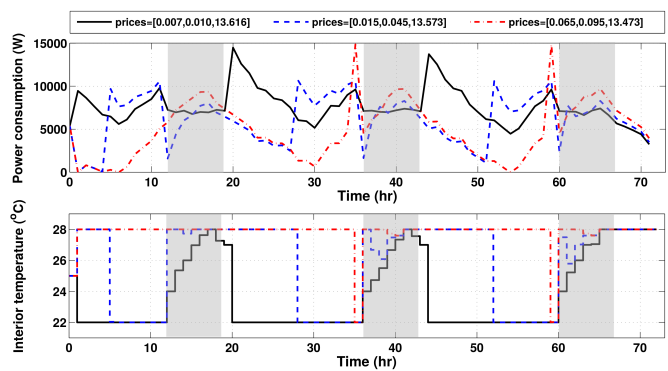}
\caption{Utility power consumption and optimal temperature settings for high, medium and low demand penalties. Shaded areas correspond to on-peak hours.}
\label{fig:scen1}
\end{figure}

To examine the impact of changes in electricity prices on peak demand, we next chose several different prices corresponding to high, medium and low penalties for peak electricity demand. Again, in each case, our algorithm is compared to GPOPS and a precooling strategy. The results are summarized in Table~\ref{tab:scen1}. For each price, the smallest computed
production cost and demand peak are typed in bold. The power consumption and the temperature settings as a function of time for the optimal strategy can be found in Fig.~\ref{fig:scen1}. For the optimal strategy, notice that by increasing the demand penalty, relative to the low-penalty case, the peak consumption is reduced by 14\% and 23\% in the medium and high penalty cases respectively. Furthermore, notice that by using the optimal strategy and the high demand-limiting prices, we have reduced the demand peak by 29\% with respect to the constant strategy in Table~\ref{tab:APS}. Of course, a moderate reduction in peak demand at the expense of large additional energy costs may not be desirable. The question of optimal distribution of electricity prices is discussed in Scenario II. \vspace{-0.1in}

\begin{table}[htbp]
\renewcommand{\tabcolsep}{4pt}
\caption{Electricity bills and demand peaks using different temperature setting strategies - Electricity prices from APS.}
\vspace{-0.1in}
\centering
\begin{tabular}{|c|c|c|c|c|}
\hline
temperature setting  & Electricity bill ($\$$) & demand peak ($kW$)  \\
\hline
\hline
Optimal (Theorem 1)& \textbf{36.58} & 9.222   \\
\hline
GPOPS~\cite{gpops}  & 37.03 & 9.155  \\
\hline
Pre-cooling        & 39.81 & \textbf{8.803}   \\
\hline
Constant           & 39.42 & 10.462    \\
\hline
\end{tabular}  \vspace{-0.1in}
\label{tab:APS}
\end{table}

\def\arraystretch{1.3}
\begin{table}[htbp]
\renewcommand{\tabcolsep}{2pt}
\caption{Demand peaks and Production costs for various prices using optimal thermostat programming, GPOPS and precooling. Marginal costs are from SRP: $[a,b]=[0.0814  \frac{\$}{kWh}, 59.76 \frac{\$}{kW}]$} \vspace{-0.1in}

\centering
\begin{tabular}{|c|c|c|c|c|c|}
\hline
& \multicolumn{1}{c|}{Prices $[p_{\text{off}}, p_{\text{on}}, p_d]$} & \multicolumn{1}{c|}{Demand-limiting} & \multicolumn{1}{c|}{Production cost} & \multicolumn{1}{c|}{Demand peak} \\
\hline
 &    $[0.007,0.010,13.616]$    &  high    & $\$$  88.712   & \textbf{7.4132} $kW$ \\
\cline{2-5}
\parbox[t]{2mm}{\multirow{-2}{*}{\rotatebox[origin=]{90}{$ \hspace{-0.2in} \text{Optimal}$}}} &    $[0.015,0.045,13.573]$    &  medium  & $\$$  \textbf{85.793}  & \textbf{8.2898} $kW$ \\
\cline{2-5}
 &    $[0.065,0.095,13.473]$    &  low     & $\$$ \textbf{86.565}   &   9.6749 $kW$ \\
\hline
\end{tabular}

\begin{tabular}{|c|c|c|c|c|c|}
\hline
\hline
 & Prices  $[p_{\text{off}}, p_{\text{on}}, p_d]$  & Demand-limiting   & Production cost  & Demand peak  \\
\hline
 &     $[0.007,0.010,13.616]$    &  high    & $\$$ \textbf{84.396}  & 7.9440 $kW$\\
\cline{2-5}
\parbox[t]{2mm}{\multirow{-2}{*}{\rotatebox[origin=]{90}{$ \hspace{-0.2in} \text{GPOPS}$}}}  &    $[0.015,0.045,13.573]$    &  medium  & $\$ $ 86.182  &  9.1486 $kW$\\
\cline{2-5}
 &     $[0.065,0.095,13.473]$    &  low     & $\$ $ 87.382   & 9.6221 $kW$\\
\hline
\end{tabular}

\begin{tabular}{|c|c|c|c|c|c|}
\hline
\hline
 &  Prices  $[p_{\text{off}}, p_{\text{on}}, p_d]$  & Demand-limiting   & Production cost  & Demand peak  \\
\hline
 &     $[0.007,0.010,13.616]$    &  high    & $ \$ $  91.064    &  8.8031 $kW$ \\
\cline{2-5}
\parbox[t]{2mm}{\multirow{-2}{*}{\rotatebox[origin=]{90}{$ \hspace{-0.15in} \text{Precooling}$}}} &     $[0.015,0.045,13.573]$    &  medium  & $ \$ $  91.064   & 8.8031 $kW$\\
\cline{2-5}
 &     $[0.065,0.095,13.473]$    &  low     & $ \$ $ 91.064  & \textbf{8.8031} $kW$\\
\hline
\end{tabular}
\label{tab:scen1} \vspace{-0.1in}
\end{table}

\subsection{Scenario 2: Optimal Thermostat Programming with Optimal Electricity Prices}
In this case, we applied Algorithm~\ref{algorithm} to find optimal on-peak, off-peak and demand prices under the assumption that the building's parameters in Table~\ref{tab:params} represent an averaged user. The marginal production costs $a$ and $b$ are taken as $[a,b] = [0.0814,59.76]$ as estimated by SRP. The optimal prices, associated production cost, and associated peak demand are listed in Table~\ref{tab:scen2}. A typical pricing strategy for SRP and other utilities is to set prices proportional to marginal production costs. The production cost associated with this strategy is also listed in Table~\ref{tab:scen2}. Notice that the optimal prices are in fact not proportional to the marginal costs of generation. \vspace{-0.1in}

\begin{table}[htbp]
\renewcommand{\tabcolsep}{2pt}
\caption{Costs of production and demand peaks associated with regulated optimal and SRP's electricity prices. Marginal costs from SRP: $a=0.0814  \frac{\$}{kWh}, b=59.76 \frac{\$}{kW}$} \vspace{-0.1in}
\centering
\begin{tabular}{|c|c|c|c|}
\hline
 Strategy &  $[p_{\text{off}} (\frac{\$}{kWh}), p_{\text{on}} (\frac{\$}{kWh}), p_d (\frac{\$}{kW})]$ &  Production cost  &  Demand peak \\
 \hline
 Optimal & $[0.082, 0.108, 54.004]$  &  $\$$ \textbf{83.333}  &  8.3008 $kW$  \\
 \hline
SRP & $[0.0572,0.0814,59.76]$  &  $\$$ 89.005  &  7.4661 $kW$ \\
\hline
\end{tabular}
\label{tab:scen2} 
\end{table}

\section{Conclusion and Future Work}
In this work, we proposed a dynamic-programming-based algorithm for solving the optimal control problem associated with thermostat programming in the presence of distributed thermal energy storage in interior structures. We used a pricing model which is a combination of on-peak, off-peak and demand charges. Using the solution to this optimal control problem as a model of behavior, we determined the optimal prices which minimize production costs for the utility. We concluded that optimal thermostat programming can significantly reduce electricity bills and demand peak by taking advantage of energy storage using thermal mass. Furthermore, we showed that the typical approach to electricity pricing is suboptimal at reducing production costs. The results of this paper assume a rational consumer and accurate models of both the daily temperature and utility production costs.

\section{Acknowledgements}
This material is based upon work supported by the National Science Foundation under Grant No. CMMI-1301851. We would like to thank Salt River Project (SRP) for providing us with their suggestions and data.

\bibliographystyle{ieeetr}
\bibliography{Reza_ACC2015}

\end{document}